\documentclass{article}
\usepackage{fullpage}
\usepackage[T1]{fontenc}
\usepackage{amsmath, amsthm, amsfonts, amssymb}
\usepackage{graphicx}
\usepackage{color}
\usepackage{lineno}
\usepackage{enumitem}
\usepackage[dvipsnames]{xcolor}
\usepackage[utf8]{inputenc}
\usepackage[english]{babel}

\newtheorem{prob}{Problem}
\newtheorem{theorem}{Theorem}[section]
\newtheorem{thm}{Theorem}
\newtheorem{cor}[thm]{Corollary}

\newtheorem*{prop}{Proposition}

\theoremstyle{definition}

\theoremstyle{remark}

\title{Hamiltonicity of the Double Vertex Graph and the \\ Complete Double Vertex Graph of some Join Graphs}
\author{Luis Adame\thanks{Unidad Acad\'emica de Matem\'aticas, Universidad Aut\'onoma de Zacatecas, Zacatecas, Mexico. \texttt{l\_e\_a\_m\_@hotmail.com, luismanuel.rivera@gmail.com}} \and Luis Manuel Rivera\footnotemark[1] \and Ana Laura Trujillo-Negrete\thanks{Departamento de Matem\'aticas, Cinvestav, CDMX, Mexico. Partially supported by CONACYT (Mexico), grant 253261. \texttt{ltrujillo@math.cinvestav.mx}}
}
\date{\today}

\begin{document}
\maketitle


\abstract{Let $G$ be a simple graph of order $n$. The double vertex graph $F_2(G)$ of $G$ is the graph whose vertices are the $2$-subsets of $V(G)$, where two vertices are adjacent in $F_2(G)$ if their symmetric difference is a pair of adjacent vertices in $G$. A generalization of this graph is the complete double vertex graph $M_2(G)$ of $G$, defined as the graph whose vertices are the $2$-multisubsets of $V(G)$, and two of such vertices are adjacent in $M_2(G)$ if their symmetric difference (as multisets) is a pair of adjacent vertices in $G$. In this paper we exhibit an infinite family of graphs (containing Hamiltonian and non-Hamiltonian graphs) for which $F_2(G)$ and $M_2(G)$ are Hamiltonian. 
This family of graphs is the set of join graphs $G=G_1 + G_2$, where $G_1$ and $G_2$ are of order $m\geq 1$ and $n\geq 2$, respectively, and $G_2$ has a Hamiltonian path. For this family of graphs, we show that if $m\leq 2n$ then $F_2(G)$ is Hamiltonian, and if $m\leq 2(n-1)$ then $M_2(G)$ is Hamiltonian. }

\section{Introduction.}

\begin{figure}[t]
	\centering
	\includegraphics[width=0.7\textwidth]{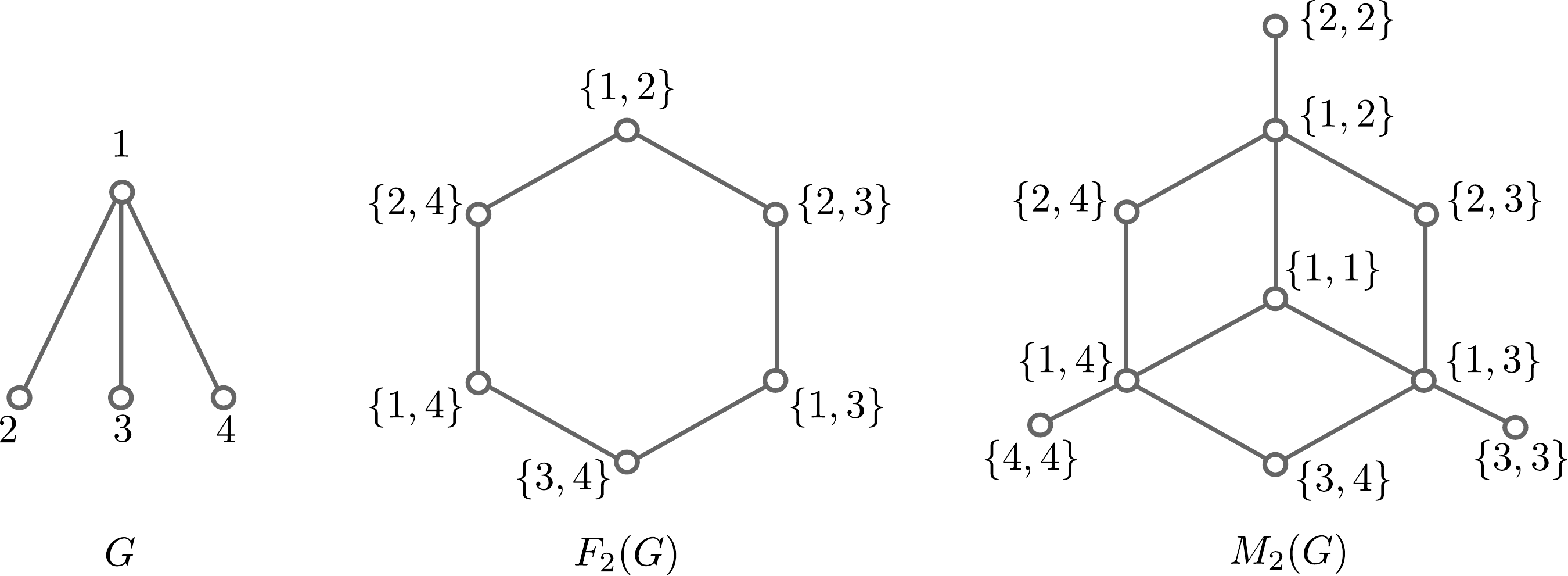}
	\caption{\small A graph $G$, its double vertex graph $F_2(G)$ and its complete double vertex graph $M_2(G)$.}
	\label{fig:definition-example}
\end{figure}

Throughout this paper, $G$ is a simple graph of order $n \geq 2$. In this paper we deal with two constructions of graphs, the 
double vertex graph and the complete double vertex graph. 
The \emph{$k$-token graph $F_k(G)$} of $G$ is the graph whose vertices are the $k$-subsets of $V(G)$, 
where two of such vertices are adjacent if their symmetric difference is a pair of adjacent vertices in $G$. 
The \emph{$k$-multiset graph $M_k(G)$} of $G$ is the graph whose vertices are the $k$-multisubsets of $V(G)$, 
and two of such vertices are adjacent if their symmetric difference (as multisets) is a pair of adjacent vertices in $G$.  
See
an example of these constructions in Figure~\ref{fig:definition-example}. 
The $2$-token graph is usually called the \emph{double vertex graph} and the $2$-multiset graph is called the \emph{complete double vertex graph}.


The $k$-token graphs have been defined, independently, at least four times, see \cite{alavi1,FFHH, johns-1988, rudolph}.  
A classical example of token graphs is the Johnson graph $J(n, k)$ that is, in fact, the $k$-token graph of 
the complete graph $K_n$. The Johnson graphs have been widely studied in the last three decades due to 
its connections with coding theory, see for example \cite{etzion, spectra-johnson, automorphisms-johnson}.   
The $k$-multiset graph was introduced in 2001 by Chartrand et al. \cite{char}. 

In 1988, Johns defined the $k$-token graphs in his PhD thesis under the name of the \emph{$k$-subgraph graph}, and he studied some combinatorial properties of these graphs. 

In 1991, Alavi et al. reintroduced, independently, the $2$-token graphs, calling them the double vertex graphs,  and they studied combinatorial properties of these graphs, such as connectivity, planarity, 
regularity and Hamiltonicity, see \cite{alavi1, alavi2, alavi3, alavi4, zhu}. 

Several years later, Rudolph \cite{aude, rudolph} redefined the token graphs, with the name of {\it symmetric powers of graphs},
with the aim to study the graph isomorphism problem and for its possible applications to quantum mechanics. Rudolph gave 
several examples of cospectral non-isomorphic graphs such that the corresponding $2$-token graphs are non-cospectral. 
This shows that, sometimes, the spectrum of the $2$-token graph of $G$ is a better invariant than the spectrum of $G$. 
However, Alzaga et al. \cite{alzaga} and, independently, Barghi and Ponomarenko \cite{barghi} proved that for any positive 
integer $k\geq 2$ there exists infinitely many pairs of non-isomorphic graphs with cospectral $k$-token graphs. 
Several authors have continued with the study of the possible applications of the token graphs in physics (see. e.g., \cite{fisch2, fisch, ouy}).

Fabila-Monroy et al.,~\cite{FFHH} reintroduced the concept of $k$-token graph of $G$ as a model in which $k$ 
indistinguishable tokens move on the vertices of a graph $G$ along the edges of $G$. They began 
a systematic study of some combinatorial parameters of $F_k(G)$ such as connectivity, diameter, cliques, chromatic number, 
Hamiltonian paths and Cartesian product. This line of research has been continued by different authors, see, e.g.,  
\cite{dealba2, deepa, deepa2, ruyanalea, soto, leancri, leatrujillo}. 
In particular Soto et al.,~\cite{soto} showed that 
a problem in coding theory is equivalent to the study of the packing number of the token graphs of the path graph.

For two disjoint graphs $G_1$ and $G_2$, the \emph{join graph}  $G=G_1 + G_2$  of graphs $G_1$ and $G_2$ is the graph 
whose vertex set is $V(G_1)\cup V(G_2)$ and its edge set is $E(G_1)\cup E(G_2)\cup \{uv:u\in G_1\text{ and }v\in G_2\}$, 
a simple example is the complete bipartite graph $K_{m,n}=E_m + E_n$, where $E_r$ denotes the graph of $r$ isolated vertices. 
The graph $F_{m,n}=E_m + P_n$ is called the \emph{fan graph}, where $P_n$ denotes the path graph of $n$ vertices (see an 
example in Figure~\ref{fig1} (left)). 


A \emph{Hamiltonian path} (resp. a \emph{Hamiltonian cycle}) 
of a graph $G$ is a path (resp. cycle) containing each vertex of $G$ exactly once. 
A graph $G$ is \emph{Hamiltonian} if it contains a Hamiltonian cycle. 
In this paper we show the following result for the fan graphs. 

\begin{thm}
	\label{thm:main}
	Let $m\geq 1$ and $n\geq 2$. Then, $F_2(F_{m,n})$ is Hamiltonian if and only if $m\leq 2n$, and  $M_2(F_{m,n})$
	is Hamiltonian if and only if $m\leq 2(n-1)$.  
\end{thm}

With the aim of clarity in the exposition of the proof, this theorem has been separated in next subsection as Theorem~\ref{thm:main1} and Theorem~\ref{thm:complete_double_vertex}. The proof of Theorem~\ref{thm:main1} was  published in~\cite{adame} (Theorem 1) and the proof of Theorem~\ref{thm:complete_double_vertex} was published in~\cite{complete-double}.

This theorem implies easily the following more general result. 
\begin{cor}
	\label{cor:main}
	Let $G_1$ and $G_2$ be two graphs
	of order $m\geq 1$ and $n\geq 2$, respectively, such that $G_2$ has a Hamiltonian path. Let $G=G_1+G_2$. 
	If $m\leq 2n$ then $F_2(G)$ is Hamiltonian, and if $m\leq 2(n-1)$ then $M_2(G)$ is Hamiltonian. 
\end{cor}

So far, the families of graphs for which it has been studied the Hamiltonicity of their double
vertex graphs and their complete double vertex graphs are the following: complete bipartite graphs or graphs that have a Hamiltonian path.  
We point out that the infinite family of graphs given by Corollary~\ref{cor:main}
contains an infinity number of non-Hamiltonian graphs for which their double vertex graphs
and complete double vertex graphs are Hamiltonian, for example, as we are going to show, if $n+1\leq m\leq2n$ (resp. $n+1\leq m\leq 2(n-1)$) 
then $F_{m,n}$ is non-Hamiltonian while its double vertex graph $F_2(F_{m,n})$ 
(resp. its complete double vertex graph $M_2(F_{m,n})$) is Hamiltonian.

\subsection{Hamiltonicity of double and complete double vertex graphs}

It is well known that the Hamiltonicity of $G$ does not imply the Hamiltonicity of $F_k(G)$. 
For example, for the complete bipartite graph $K_{m,m}$, Fabila-Monroy et al.~\cite{FFHH} showed that if $k$ is even, then 
$F_k(K_{m, m})$ is non-Hamiltonian. 
A more easy and traditional example is the case of a cycle graph. 
It is known that if $n=4$ or $n \geq 6$, then $F_2(C_n)$ is not Hamiltonian. On the other hand, 
there exist non-Hamiltonian graphs for which its double vertex graph is Hamiltonian,
a simple example is the graph $K_{1,3}$, for which $F_2(K_{1,3})\simeq C_6$, and so 
$F_2(K_{1,3})$ is Hamiltonian. 

Next, we list the known results about the Hamiltonicity of $F_k(G)$ or the existence of a Hamiltonian path in $F_k(G)$, 
when $k$ may be greater than two. 
\begin{itemize}
	\item If $n\geq 3$ and $1\leq k\leq n-1$, then $F_k(K_n)$ is Hamiltonian, see for example \cite{alspach}.
	\item If $m\geq2$, then $F_k(K_{m,m})$ has a Hamiltonian path if and only if $k$ is odd \cite{FFHH}. 
	\item If $G$ is a graph containing a Hamiltonian path and $n$ is even and $k$ is odd, then $F_k(G)$ has a Hamiltonian path \cite{FFHH}. 
	\item If  $n\geq 3$ and $1\leq k\leq n-1$, then $F_k(F_{1,n-1})$ is Hamiltonian \cite{rive-tru}.   
\end{itemize}
In addition to these results, the following are some known results for the double vertex graph ($k=2$). 
\begin{itemize}
	\item $F_2(C_n)$ is non-Hamiltonian \cite{alavi4}. 
	\item If $G$ is a cycle with an odd chord, then $F_2(G)$ is Hamiltonian \cite{alavi4}. 
	\item $F_2(K_{m,n})$ is Hamiltonian if and only if $(m-n)^2=m+n$ \cite{alavi4}. 
\end{itemize}
More results about the Hamiltonicity of double vertex graphs can be found in the survey of Alavi et. al.~\cite{alavi3}.
We point out that the graphs for which have been studied the Hamiltonicity of its $k$-token graph (even for $k=2$)
are Hamiltonian or have a Hamiltonian path or are complete bipartite graphs.

As we mentioned before, in 2018 \cite{rive-tru} the last two authors of this article showed the following result: 
if $n \geq 3$, and $1\leq k \leq n-1$, then the $k$-token graph of the fan graph $F_{1, n-1}$ 
is Hamiltonian. We have continued with this line of research and 
in this work we show the following result for the double vertex graph of fan graphs.

\begin{theorem}[\cite{adame}, Theorem 1]
	\label{thm:main1}
	The double vertex graph of $F_{m,n}$ is Hamiltonian if and only if $n\geq 2$ and $1\leq m \leq 2\,n$, or $n=1$ and $m=3$. 
\end{theorem}

\begin{figure}[h]
\begin{center}
\includegraphics[angle=0, width=14cm]{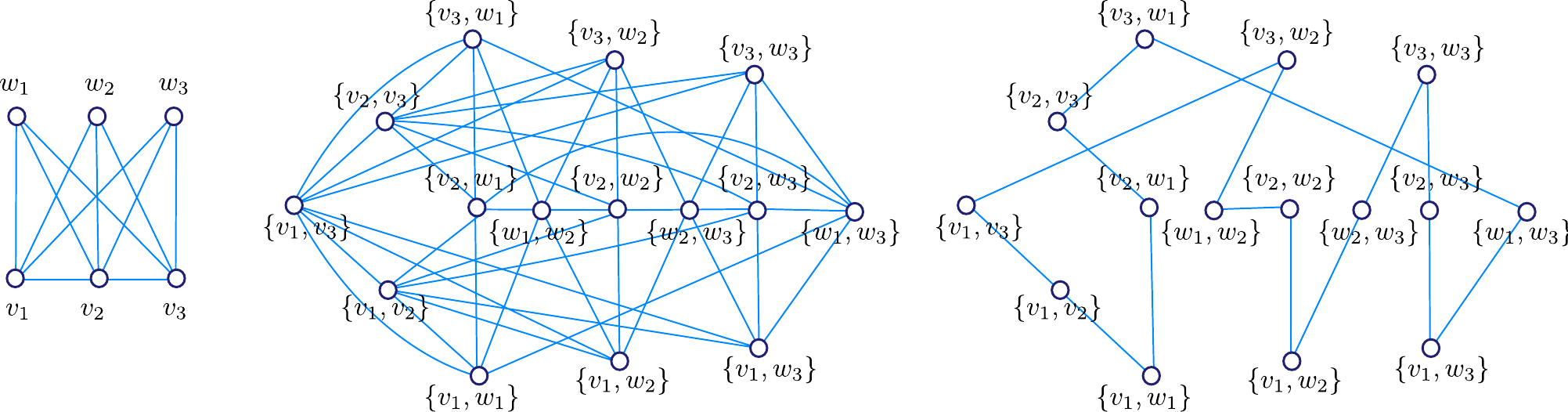}
\caption{\small The graph $F_{3,3}$ (left), its double vertex graph (center) and a Hamiltonian cycle in such graph (right).}
\label{fig1}
\end{center}
\end{figure}

In Figure~\ref{fig1} we show the double vertex graph of the fan graph $F_{3, 3}$ (center) and a Hamiltonian cycle in such graph (right). 

Let us now turn our attention to the complete double vertex graph. 
Complete double vertex graphs were implicitly presented in the work of  Chartrand et al. \cite{char}, 
and in an explicit way by Jacob et. al. \cite{jacob}, were their  first combinatorial properties were studied, 
and are a generalization of the double vertex graphs. 

As far as we know, only the following two results are known about the Hamiltonicity of $M_2(G)$. 
\begin{itemize}
	\item For $n\geq 4$, $M_2(C_n)$ is non-Hamiltonian \cite{jacob}. 
	\item If $G$ is obtained from a cycle of $n$ vertices by adding a chord between two vertices at distance two, 
	then $M_2(G)$ is Hamiltonian \cite{jacob}. 
\end{itemize}


In this paper we show the following  result for the Hamiltonicity of complete double vertex graphs.

\begin{theorem}[\cite{complete-double}, Theorem 1]
	\label{thm:complete_double_vertex}
	The complete double vertex graph of $F_{m,n}$ is Hamiltonian if and only if $n\geq 2$ and $1\leq m \leq 2\,(n-1)$. 
\end{theorem}

As we mentioned before, we have splitted Theorem~\ref{thm:main} into Theorems~\ref{thm:main1} and \ref{thm:complete_double_vertex}, 
and Corollary~\ref{cor:main} follows easily from Theorem~\ref{thm:main}. The infinite family of graphs
given in Corollary~\ref{cor:main} contains an infinite number of non-Hamiltonian graphs for which their double
vertex graph and complete double vertex graph are Hamiltonian, for example, for the fan graph $F_{m,n}$ we know
that $F_{m,n}$ is Hamiltonian if and only if $1\leq m\leq n$, while, as we are going to show, $F_2(F_{m,n})$ (resp. $M_2(F_{m,n})$)
is Hamiltonian if and only if $1\leq m\leq 2n$ (resp. $1\leq m\leq 2(n-1)$). 

%
%
%

The rest of the paper is organized as follows. In Section~\ref{sec:double} we
present the proof of Theorem~\ref{thm:main1} and in Section~\ref{sec:complete} the proof of Theorem~\ref{thm:complete_double_vertex};
our strategy to prove these results is to show explicit Hamiltonian cycles in each case.
For the purpose of clarity, in Section~\ref{sec:examples} we present some examples of our constructions.  
Finally, we suggest some open problems in Section~\ref{sec:open}. 

Before go further, let us establish some notation. Let $V(P_n):=\{v_1,\ldots,v_n\}$ 
and $V(E_m):=\{w_1,\ldots, w_m\}$, so we have  
$V(F_{m,n})=\{v_1,\ldots,v_n,w_1,\ldots,w_m\}$. 
For a path $T=a_1 a_2 \dots a_{l-1}a_{l}$, we denote by $\overleftarrow{T}$ to the reverse path $a_la_{l-1}\dots a_2a_1$. 
As usual, for a positive integer $r$, we denote by $[r]$ to the set $\{1,2,\ldots,r\}$. 
For a graph $G$, we denote by $\mu(G)$ to the number of components of $G$.

\section{Proof of Theorem~\ref{thm:main1}}
\label{sec:double}


If $n=1$ then $G\simeq K_{1,m}$, and it is know that $F_2(K_{1,m})$ is Hamiltonian if and only if $m=3$ 
(see, e.g., Proposition 5 in~\cite{alavi3}). From now on, assume $n\geq 2$.  We distinguish four cases: either $m=1$, $m=2n$, 
	$1<m<2n$ or $m>2n$. 
	
	\begin{itemize}
		\item {\bf Case} $\mathbf{m=1.}$ 	
			
		For $n=2$ we have $F_2(F_{1,2})\simeq F_{1,2}$, and so $F_2(F_{1,2})$ is Hamiltonian. 
		Now we work the case $n\geq 3$. For $1\leq i< n$ let 
		\[T_i:=\{v_i,w_1\}\{v_i,v_{i+1}\}\{v_i, v_{i+2}\} \dots \{v_i,v_n\}\]
	    and let
	   \[
	   T_n:=\{v_n,w_1\}.
	   \]
		It is clear that every $T_i$ is a path in $F_2(F_{1,n})$ and that $\{T_1, \dots, T_n\}$ 
		is a partition of $V(F_2(F_{1,n}))$. 
		
Let 
			\[
			C:=
			\begin{cases}
			\overleftarrow{T_1}\,T_2\overleftarrow{T_3}\,T_4\dots \overleftarrow{T_{n-1}}\,T_n(v_1, v_n) & \text{if $n$ is even,} \\
			\overleftarrow{T_1}\,T_2\overleftarrow{T_3}\,T_4\,\dots \,T_{n-1}\,\overleftarrow{T_{n}}\,\{v_1, v_n\} & \text{if $n$ is odd.}
			\end{cases}
			\]
			
			We are going to show that $C$ is a Hamiltonian cycle of $F_2(F_{1,n})$. Suppose $n$ is even, so 
			\[
			\begin{aligned}
			\begin{split}
			C=&\underbrace{\{v_1, v_n\} \dots \{v_1, w_1\}}_{\overleftarrow{T_1}}\underbrace{\{v_2, w_1\} \dots \{v_2, v_n\}}_{T_2}\underbrace{\{v_3, v_n\} \dots \{v_3, w_1\}}_{\overleftarrow{T_3}}\;\dots \\
			&	\underbrace{\{v_{n-1}, v_n\}\{v_{n-1}, w_1\}}_{\overleftarrow{T_{n-1}}}\underbrace{\{v_n, w_1\}}_{T_n}\{v_1, v_n\}.
			\end{split}
			\end{aligned}\]

			We are going to show that $C$ is a Hamiltonian cycle of $F_2(F_{1,n})$. First, note that for $i$ odd, the final
			vertex of $\overleftarrow{T_i}$ is $\{v_i,w_1\}$, while the initial vertex of $T_{i+1}$ is $\{v_{i+1},w_1\}$, and
			since these two vertices are adjacent in $F_2(F_{1,n})$, the concatenation $\overleftarrow{T_i}\,T_{i+1}$ corresponds to 
			a path in $F_2(F_{1,n})$. Similarly, for $1\leq i<n$ even, the final vertex of $T_i$ is $\{v_i,v_n\}$ while the initial
			vertex of $\overleftarrow{T_{i+1}}$ is $\{v_{i+1},v_n\}$, so again, the concatenation $T_i\,\overleftarrow{T_{i+1}}$
			corresponds to a path in $F_2(F_{1,n})$. Also note that the unique  vertex of $T_n$ is $\{v_n,w_1\}$ that  is adjacent, in $F_2(F_{1,n})$, to $\{v_1,v_n\}$. As the first vertex of  $\overleftarrow{T_1}$ is $\{v_1,v_n\}$, we have that $C$ is a cycle in $F_2(F_{1,n})$. 
			
			
			 Case $n$ odd.

			That is
			\[
			\begin{aligned}
			\begin{split}
			C=&\underbrace{\{v_1, v_n\} \dots \{v_1, w_1\}}_{\overleftarrow{T_1}}\underbrace{\{v_2, w_1\} \dots \{v_2, v_n\}}_{T_2}\underbrace{\{v_3, v_n\} \dots \{v_3, w_1\}}_{\overleftarrow{T_3}}\dots \\
			&	\underbrace{\{v_{n-1}, w_1\}\{v_{n-1}, v_n\}}_{T_{n-1}}\underbrace{\{v_n, w_1\}}_{\overleftarrow{T_n}}\{v_1, v_n\}.
			\end{split}
			\end{aligned}\]	
			
			In a similar way to the previous case, we can prove that $C$ is a Hamiltonian cycle of $F_2(F_{1,n})$.


Now we construct a Hamiltonian cycle in $F_2(F_{1,n})$ depending on   the parity of $n$. 

		\begin{itemize}
			\item Case $n$ even. 
			
			Let 
			\[
			C:=\overleftarrow{T_1}\,T_2\overleftarrow{T_3}\,T_4\dots \overleftarrow{T_{n-1}}\,T_n(v_1, v_n). 
			\]
			
That is 
\[
\begin{aligned}
\begin{split}
C=&\underbrace{\{v_1, v_n\} \dots \{v_1, w_1\}}_{\overleftarrow{T_1}}\underbrace{\{v_2, w_1\} \dots \{v_2, v_n\}}_{T_2}\underbrace{\{v_3, v_n\} \dots \{v_3, w_1\}}_{\overleftarrow{T_3}}\;\dots \\
&	\underbrace{\{v_{n-1}, v_n\}\{v_{n-1}, w_1\}}_{\overleftarrow{T_{n-1}}}\underbrace{\{v_n, w_1\}}_{T_n}\{v_1, v_n\}.
\end{split}
\end{aligned}\]

			We are going to show that $C$ is a Hamiltonian cycle of $F_2(F_{1,n})$. First, note that for $i$ odd, the final
			vertex of $\overleftarrow{T_i}$ is $\{v_i,w_1\}$, while the initial vertex of $T_{i+1}$ is $\{v_{i+1},w_1\}$, and
			since these two vertices are adjacent in $F_2(F_{1,n})$, the concatenation $\overleftarrow{T_i}\,T_{i+1}$ corresponds to 
			a path in $F_2(F_{1,n})$. Similarly, for $1\leq i<n$ even, the final vertex of $T_i$ is $\{v_i,v_n\}$ while the initial
			vertex of $\overleftarrow{T_{i+1}}$ is $\{v_{i+1},v_n\}$, so again, the concatenation $T_i\,\overleftarrow{T_{i+1}}$
			corresponds to a path in $F_2(F_{1,n})$. Also note that the unique  vertex of $T_n$ is $\{v_n,w_1\}$ that  is adjacent, in $F_2(F_{1,n})$, to $\{v_1,v_n\}$. As the first vertex of  $\overleftarrow{T_1}$ is $\{v_1,v_n\}$, we have that $C$ is a cycle in $F_2(F_{1,n})$. 
			
			
			\item  Case $n$ odd. 
			
			Let 			
			\[C:=\overleftarrow{T_1}\,T_2\overleftarrow{T_3}\,T_4\,\dots \,T_{n-1}\,\overleftarrow{T_{n}}\,\{v_1, v_n\}.\]
			
	That is
\[
\begin{aligned}
\begin{split}
	C=&\underbrace{\{v_1, v_n\} \dots \{v_1, w_1\}}_{\overleftarrow{T_1}}\underbrace{\{v_2, w_1\} \dots \{v_2, v_n\}}_{T_2}\underbrace{\{v_3, v_n\} \dots \{v_3, w_1\}}_{\overleftarrow{T_3}}\dots \\
&	\underbrace{\{v_{n-1}, w_1\}\{v_{n-1}, v_n\}}_{T_{n-1}}\underbrace{\{v_n, w_1\}}_{\overleftarrow{T_n}}\{v_1, v_n\}.
\end{split}
\end{aligned}\]	

In a similar way to the previous case, we can prove that $C$ is a Hamiltonian cycle of $F_2(F_{1,n})$. 
\end{itemize}

\item {\bf Case} $\mathbf{m=2n.}$

 Let $C$ be the cycle defined in the previous case depending on the parity of $n$.
Let 
 \[
 P_1:=\{v_n,w_1\} \xrightarrow{C} \{v_1,v_n\}
 \]   
 be the path from $\{v_n, w_1\}$ to $\{v_1, v_n\}$ 
obtained from $C$ by deleting the edge between $\{v_n, w_1\}$ and
$\{v_1, v_n\}$. That is 
\[P_1=\begin{cases}
\{v_n,w_1\}\,T_{n-1}\,\dots\,\overleftarrow{T_4}\,T_3\,\overleftarrow{T_2}\,T_1 & \text{if $n$ is even, }\\
\{v_n, w_1\}\,\overleftarrow{T_{n-1}}\, \dots\,\overleftarrow{T_4}\,T_3\,\overleftarrow{T_2}\,T_1 & \text{if $n$ is odd. }\\
\end{cases}
\]

For $1< i\leq n$ let
\[
\begin{split}
P_i:=&\{w_i,v_n\}\{w_i,w_1\}\{w_i,v_{n-1}\}\{w_i,w_{i+(n-1)}\}\{w_i,v_{n-2}\}\{w_i,w_{i+(n-2)}\}\{w_i,v_{n-3}\}\{w_i,w_{i+(n-3)}\}\ldots \\
&\{w_i,v_2\}\{w_i,w_{i+2}\}\{w_i,v_1\}\{w_i,w_{i+1}\}.
\end{split}
\]
We can observe that after $\{w_i, w_1\}$ the vertices in $P_i$ follows the  pattern $\{w_{i}, v_j\}\{w_{i}, w_{i+j}\}$, from $j=n-1$ to $1$.
For $n+1\leq i\leq 2n$ let
\[
\begin{split}
P_i:=&\{w_i,v_n\}\{w_i,w_{i+n}\}\{w_i,v_{n-1}\}\{w_i,w_{i+(n-1)}\}\{w_i,v_{n-2}\}\{w_i,w_{i+(n-2)}\}\ldots\\
& \{w_i,v_2\}\{w_i,w_{i+2}\}\{w_i,v_1\}\{w_i,w_{i+1}\},
\end{split}
\]
where the sums are taken mod $2n$ with the convention that $2n\pmod{2n}=2n$. In this case, the vertices  in $P_i$ after $\{w_i,w_{i+n}\}$ follows the  pattern $\{w_{i}, v_j\}\{w_{i}, w_{i+j}\}$, from $j=n-1$ to $1$.

We claim that the concatenation 
\[
C_2:=P_1\,P_2\,\ldots\,P_{2n}\{v_n,w_1\}
\]
is a Hamiltonian cycle in $F_2(F_{m,n})$. First we prove that $\{P_1,\ldots,P_{2n}\}$ is a partition of $F_2(F_{m,n})$. It is clear that the paths $P_1,\ldots,P_{2n}$ are pairwise disjoint in $F_2(F_{m,n})$. Now, we are going to show that every vertex in $F_2(F_{m,n})$ belongs to exactly one of the paths $P_1,\ldots,P_{2n}$.  

\begin{itemize}
	\item $\{v_i, v_j\}$ belongs to $P_1$, for any $i,j\in [n]$ with $i\neq j$.
	\item $\{w_i,v_j\}$ belongs to $P_i$, for any $i\in [m]$ and $j\in [n]$.
	\item $\{w_i, w_1\}$ belong to $P_i$, for any $i\in [m]$. 
	\item Consider now the vertices of type $\{w_i,w_j\}$, for $1<i<j\leq n$, 
	\begin{itemize}
		\item $\{w_i,w_j\}$ belongs to $P_i$, for any $1<i\leq n$ and $i<j\leq i+n-1$.
		\item $\{w_i,w_j\}$ belongs to $P_j$, for any $1<i\leq n$ and $i+n-1< j\leq 2n$.
		\item $\{w_i,w_j\}$ belongs to $P_i$, for any $n< i< 2n$ and $i<j\leq 2n$.
	\end{itemize}
\end{itemize}

Now we show that $C$ is a cycle in $F_2(F_{m,n})$. Note that
the final vertex of $P_1$ is $\{v_1,v_n\}$ while the initial vertex of $P_2$ is $\{w_2, v_n\}$, and these
two vertices are adjacent in $F_2(F_{m,n})$. Also, for $1<i<2n$, the final vertex of $P_i$ is $\{w_i,w_{i+1}\}$ while
the initial vertex of $P_{i+1}$ is $\{w_{i+1},v_n\}$, and again these two vertices are adjacent in $F_2(F_{m,n})$. 
On the other hand, the final vertex of $P_{2n}$ is $\{w_1,w_n\}$ while the initial vertex of $P_1$ is $\{v_n,w_1\}$, and
these two vertices are adjacent in $F_2(F_{m,n})$. These four observations together imply that $C_2$ is a cycle in $F_2(F_{m,n})$, and hence, $C_2$ is a Hamiltonian cycle of $F_2(F_{m,n})$. 

\item {\bf Case} $\mathbf{1<m<2\,n.}$

Consider again the paths $P_1,\dots, P_m$ defined in the previous case and let us modify them slightly in the following
way:  
\begin{itemize}
	\item $P_1'=P_1$;
	\item for $1<i<m$, let $P'_i$ be the path obtained from $P_i$ by deleting the vertices of type $\{w_i, w_j\}$, for each $j>m$;
	\item let $P_m'$ be the path obtained from $P_m$ by first interchanging the vertices $\{w_m, w_{m+1}\}$ and $\{w_m, w_1\}$ from their current positions in $P_m$, and then deleting the vertices of type $\{w_m, w_j\}$, for every $j>m$.
\end{itemize}
Given this construction of $P'_i$ we have the following: 
\begin{itemize}
	\item[(A1)] $P'_i$ induces a path in $F_2(F_{m,n})$;
	\item[(A2)] for $1\leq i<m$ the path $P'_i$ has the same initial and final vertices as the path $P_i$, and $P'_m$ has the same initial vertex as $P_m$, and its final vertex is $\{w_i, w_1\}$; 
	\item[(A3)] since we have deleted only the vertices of type $\{w_i, w_j\}$ from $P_i$ to obtain $P_i'$, for each $j>m$ and $i \in [m]$, it follows that $\{V(P'_1),\ldots,V(P'_m)\}$ is a partition of $V(F_2(F_{m,n}))$. 
	
\end{itemize}

By (A1) and (A2) we can concatenate the paths $P'_1,\ldots,P'_m$  into a cycle $C'$ as follows:
\[
C':=P'_1\,P'_2\,\ldots\,P'_m(v_n, w_1)
\]
and then by (A3) it follows that $C'$ is a Hamiltonian cycle in $F_2(F_{m,n})$.

\item {\bf Case} $\mathbf{m>2\,n.}$

Here, our aim is to show that $F_2(F_{m,n})$ is not Hamiltonian by using 
the following known result posed in West's book~\cite{west}. 

\begin{prop}[Prop. 7.2.3, \cite{west}] If $G$ has a Hamiltonian cycle, then for each 
	nonempty set $S\subset V(G)$, the graph $G-S$ has at most $|S|$ components.  
\end{prop}

Then, we are going to exhibit a subset $A\subset V(F_2(F_{m,n}))$ such that
\[
\mu(F_2(F_{m,n})-A)>|A|.
\]

Let 
\[
A:=\big\{\{w_i,v_j\}:i\in [m] \text{ and }j\in [n]\big\}.
\]

Note that for any $i,j\in [m]$ with $i\neq j$, $\{w_i,w_j\}$ is an isolated vertex of $F_2(F_{m,n})-A$, 
and there are $\binom{m}{2}$ vertices of this type. Also note that
the subgraph induced by the vertices of type $\{v_i,v_j\}$, for $i,j\in [n]$ and $i\neq j$, 
is a component of $F_2(F_{m,n})-A$, and since
$|A|=mn$ and $m>2n$, we have
\[\mu(F_2(F_{m,n})-A)\geq \binom{m}{2}+1=\frac{m(m-1)}{2}+1\geq mn+1>mn=|A|,\]
as required. This completes the proof of Theorem~\ref{thm:main1}. 
\end{itemize}

\section{Proof of Theorem~\ref{thm:complete_double_vertex}}
\label{sec:complete}

If $n=1$ and $m\geq 1$ then $\deg(\{w_i,w_i\})=1$, for any $i\in [m]$, which implies that 
$M_2(F_{m,n})$ is not Hamiltonian, so we assume that $n\geq 2$.  

The constructions that we give in this section are similar to those given in the previous section.  
	We distinguish four cases: either $m=1$, $m=2\,(n-1)$, $1<m<2\,(n-1)$ or $m>2(n-1)$. 
		\begin{itemize}
		\item {\bf Case} $\mathbf{m=1.}$ 
		
		For $1\leq i\leq n$ let 
	\[
	T'_i:=\{v_i,w_1\}\{v_i,v_i\}\{v_i,v_{i+1}\}\{v_i, v_{i+2}\}\ldots \{v_i,v_n\}.
	\] 
	We remark the following: 
	\begin{itemize}
		\item[(i)] $T'_i$ can be seen as the resulting path from $T_i$ (defined in Section 2) by adding the vertex $\{v_i,v_i\}$
		between the vertices $\{v_i,w_1\}$ and $\{v_i,v_{i+1}\}$. 
		\item[(ii)] $T'_i$ is a path in $M_2(F_{m,n})$ and $T'_n=\{v_n,w_1\}\{v_n,v_n\}$.
		\item[(iii)] For $1\leq i<n$, the paths $T'_i$ and $T_i$ have the same initial and final vertices. 
		\item[(iv)] The set $\{T'_1, \dots, T'_n, \{\{w_1, w_1\}\}\}$ is a partition of $V(M_2(F_{m,n}))$. 
	\end{itemize}
	Let 
	\[
	C:= \begin{cases}
	T'_1\, \overleftarrow{T'_2}\, T'_3 \, \overleftarrow{T'_4}\,\dots\, T'_{n-1}\,\overleftarrow{T'_n} \,\{w_1,w_1\} & \text{if $n$ is even, } \\
	\overleftarrow{T'_1}\,\{w_1,w_1\}\,T'_2\,\overleftarrow{T'_3}\,T'_4\ldots\,T'_{n-1}\,\overleftarrow{T'_n} & \text{if $n$ is odd.}
	\end{cases}
	\]
	We claim that $C$ is a Hamiltonian cycle in $M_2(F_{m,n})$. Suppose that $n$ is even. Since $(ii)$ and $(iii)$ hold, we can concatenate 
	the paths $T'_1\, \overleftarrow{T'_2}\, T'_3 \, \overleftarrow{T'_4}\,\dots\, T'_{n-1}\,\overleftarrow{T'_n}$, and since 
	$\{v_n,w_1\}$ (the final vertex of $\overleftarrow{T'_n}$) is adjacent to $\{w_1,w_1\}$, and $\{w_1,w_1\}$ is adjacent to $\{v_1,v_1\}$
	(the initial vertex of $T'_1$), it follows that $C$ is a cycle in $M_2(F_{m,n})$. By similar arguments, in the case $n$ odd
	we have that $C$ is a cycle in $M_2(F_{m,n})$. Finally, in both cases, $(iv)$ implies that $C$ is a Hamiltonian
	cycle in $M_2(F_{m,n})$, as claimed. \\
Note that in both cases, the vertices $\{v_n,w_1\}$ and $\{v_n,v_n\}$ are adjacent in $C$ (these two vertices 
correspond to the vertices of $\overleftarrow{T'_n}$).  This observation will be useful in the following two cases.

\item {\bf Case} $\mathbf{m=2\,(n-1).}$
	
	Let $C$ be the cycle defined in the previous case,  depending on the parity of $n$. Let 
	\[
	P_1:=\{v_n, w_1\} \xrightarrow{C} \{v_n, v_n\}
	\]
	be the path obtained from $C$ by deleting the edge between $\{v_n, w_1\}$ and $\{v_n, v_n\}$. That is 
	\[
	P_1:=\begin{cases}
	\{v_n, w_1\}\,\{w_1, w_1\}\,T'_1\,\overleftarrow{T'_2}\,T'_3\, \dots \,T'_{n-1}\,\{v_n, v_n\} & \text{if $n$ is even,}\\
	\{v_n, w_1\}\,\overleftarrow{T'_1}\,\{w_1, w_1\}\,T'_2\,\overleftarrow{T'_3}\,\dots\,T'_{n-1}\,\{v_n, v_n\} & \text{if $n$ is odd. }
	\end{cases} 
	\]

For $1< i\leq n-1$ let
	\[
	\begin{split}
	P_i:=&\{w_i,v_n\}\{w_i,w_i\}\{w_i,v_{n-1}\}\{w_i,w_1\}\{w_i,v_{n-2}\}\{w_i,w_{i+(n-2)}\}\{w_i,v_{n-3}\}\{w_i,w_{i+(n-3)}\}\ldots\\
	&\{w_i,v_2\}\{w_i,w_{i+2}\}\{w_i,v_1\}\{w_i,w_{i+1}\}. 
	\end{split}
	\]
	Note that after $\{w_i, w_1\}$, the vertices in $P_i$ follows the  pattern $\{w_{i}, v_j\}\{w_{i}, w_{i+j}\}$, from $j=n-1$ to $1$.
For $n\leq i\leq m$ let
	\[
	\begin{split}
	P_i:=&\{w_i,v_n\}\{w_i,w_i\}\{w_i,v_{n-1}\}\{w_i,w_{i+(n-1)}\}\{w_i,v_{n-2}\}\{w_i,w_{i+(n-2)}\}\ldots\\
	& \{w_i,v_2\}\{w_i,w_{i+2}\}\{w_i, v_1\}\{w_i,w_{i+1}\},
	\end{split}
	\]
	where the sums are taken mod $m$ with the convention that $m\pmod{m}=m$. In this case, after $\{w_i, w_i\}$, the vertices in $P_i$ follows the pattern $\{w_{i}, v_j\}\{w_{i}, w_{i+j}\}$, from $j=n-1$ to $1$.
	We claim that the concatenation 
\[
P:=P_1\,P_2\,\ldots\,P_m\{v_{n}, w_{1}\}
\]
	is a Hamiltonian cycle in $M_2(F_{m,n})$. First note that  the final vertex of $P_1$ is $\{v_n,v_n\}$ while the initial vertex of $P_2$ is $\{w_2,v_n\}$, and 
	these two vertices are adjacent in $M_2(F_{m,n})$. Moreover, for $1<i<2(n-1)$, the final vertex
	of $P_i$ is $\{w_i,w_{i+1}\}$ while the initial vertex of $P_{i+1}$ is $\{w_{i+1},v_n\}$, and also
	these two vertices are adjacent in $M_2(F_{m,n})$. Also, the final vertex of $P_m$ is $\{w_m,w_1\}$ while
	the initial vertex of $P_1$ is $\{v_n,w_1\}$, and these two vertices are adjacent in  $M_2(F_{m,n})$. 
	These three observations together imply that $P$ is a cycle in $M_2(F_{m,n})$. 
	
	It remains to show that the cycle $P$ is Hamiltonian. 
	Notice that any vertex in $V(F_{m,n})$ belong to exactly one of the following options:  
	\begin{itemize}
		\item The vertices of type $\{v_i,v_j\}$ belongs to $P_1$	for any $i,j\in [n]$.
		\item The vertices of type $\{w_i,v_j\}$ belongs to $P_i$ for any $i\in [m]$ and $j\in [n]$. 
		\item The vertices of type $\{w_i,w_i\}$ and $\{w_i,w_1\}$ belong to $P_i$ for any $i\in [m]$. 
		\item Consider now the vertices of type $\{w_i,w_j\}$ for $i\neq j$, assuming without loss of generality
		that $i<j$.  
		\begin{itemize}
			\item $\{w_i,w_j\}$ belongs to $P_i$ for any $1<i<n$ and $i<j<i+n-1$.
			\item $\{w_i,w_j\}$ belongs to $P_j$ for any $1<i<n$ and $i+n-1\leq j\leq 2(n-1)$.
			\item $\{w_i,w_j\}$ belongs to $P_i$ for any $n\leq i< 2(n-1)$ and $i<j\leq 2(n-1)$.
		\end{itemize}
	\end{itemize}
    Thus, $P$ is our desired Hamiltonian cycle in $M_2(F_{m,n})$. 

\item {\bf Case}  $\mathbf{1<m<2\,(n-1).}$

We consider again the paths $P_1,\ldots,P_m$ defined in the previous case 
with a slight modification: 
\begin{itemize}
\item $P_1'=P_1$;
\item for $i \in \{2, \dots, m-1\}$, let $P'_i$ be the path obtained from $P_i$ by deleting the vertices of type $\{w_i, w_j\}$, for each $j>m$;
\item let $P'_m$ be the path obtained from $P_m$ by first interchanging $\{w_m, w_{m+1}\}$ and $\{w_m, w_1\}$ from their current positions in $P_m$, and then deleting the vertices of type $\{w_m, w_j\}$, for every $j>m$.
\end{itemize}
We have  that $P'_1,\ldots, P'_m$ are, indeed, disjoint paths in $M_2(F_{m,n})$, and that
$P'_i$ has the same initial and final vertices as $P_i$, so the concatenation 
\[
P':=P'_1\,\ldots\,P'_m\{v_n, w_1\}
\]
correspond to a cycle in $M_2(F_{m,n})$. It is an easy exercise 
(similar as in the case of double vertex graphs) to show that this cycle is in fact a Hamiltonian cycle in $M_2(F_{m,n})$. 

\item {\bf Case}  $\mathbf{m>2\,(n-1).}$

We are going to show that, in this case, $M_2(F_{m,n})$ is not Hamiltonian. 
For this, we proceed similarly to the case $m>2n$ of Section~\ref{sec:double}, so, 
we make use of  Proposition 7.2.3 posed in West's book~\cite{west}. Thus,  we are going 
to exhibit a subset $A\subset V(M_2(F_{m,n}))$
such that 
\[
\mu(M_2(F_{m,n})-A)>|A|.
\]
Let 
\[
A:=\{\{w_i,v_j\}\in M_2(F_{m,n}) : i\in [m]\text{ and }j\in [n]\},\]
\[
T:=\{\{w_i,w_j\}\in M_2(F_{m,n}) : i,j\in [m]\}
\]
 and 
\[
R:=\{\{v_i,v_j\}\in M_2(F_{m,n}) : i,j\in [n]\}.
\]
 The set $\{A, T, R\}$ is a partition of $M_2(F_{m,n})$. 
 Note that any vertex in $T$ has its neighbors in $A$, so the subgraph induced by 
 $T$ in $M_2(F_{m,n})-A$ is the empty graph of order $\binom{m+1}{2}$ ($\overline{K_{\binom{m+1}{2}}}$). 
 On the other hand, note that the subgraph of $M_2(F_{m,n})$ induced by $R$ is isomorphic to the complete double vertex
 graph of the path of $n$ vertices (which is connected), also note that the vertices in $R$ have neighbours
 in $A$ but not in $T$, implying that the subgraph induced by $R$ is a component of $M_2(F_{m,n})-A$.  
 Since $|A|=m\,n$, $|T|=\binom{m+1}{2}$ and $m>2\,(n-1)$, we have 
\[
\mu(M_2(F_{m,n})-A)= |T|+1=\binom{m+1}{2}+1> m\,n=|A|.
\]
This completes the proof of Theorem~\ref{thm:complete_double_vertex}. 
 	\end{itemize}

\section{Some examples}
\label{sec:examples}

For the purpose of clarity, we exhibit several examples of our constructions. 
\\

\noindent\textbf{Double vertex graph of $\mathbf{F_{m,n}}$}

Following the proof of Theorem~\ref{thm:main1}, we examine first the case $m=1$, then the case $m=2n$ and finally the case
$1<m<2n$. 

\begin{itemize}
	
	\item[1)] $\mathbf{m=1.}$

	In this case we consider the fan graph $F_{1,4}$, so the corresponding paths in $F_2(F_{1,4})$ are the following:
	\[
	\begin{aligned}
	T_1=&\{v_1,w_1\}\{v_1,v_2\}\{v_1,v_3\}\{v_1,v_4\} \\
	T_2=&\{v_2,w_1\}\{v_2,v_3\}\{v_2,v_4\}\\
	T_3=&\{v_3,w_1\}\{v_3,v_4\}\\
	T_4=&\{v_4,w_1\}\\
	\end{aligned}
	\]
	
	Thus, the concatenation 
	\[
	C=\underbrace{\{v_1,v_4\}\{v_1,v_3\}\{v_1,v_2\}\{v_1,w_1\}}_{\overleftarrow{T_1}}
	\underbrace{\{v_2,w_1\}\{v_2,v_3\}\{v_2,v_4\}}_{T_2}\underbrace{\{v_3,v_4\}\{v_3,w_1\}}_{\overleftarrow{T_3}}
	\underbrace{\{v_4,w_1\}}_{T_4} 
	\]
	 
	is our desired Hamiltonian cycle in  $F_2(F_{1,4})$.

	\item[2)] $\mathbf{m=2n.}$ 
	
	Here we consider the graph $F_{8,4}$, so in $F_2(F_{8,4})$ we have the paths: 
	\[
	\begin{aligned}[]
	P_1=&\{v_4,w_1\}\{v_3,w_1\}\{v_3,v_4\}\{v_2,v_4\}\{v_2,v_3\}\{v_2,w_1\}\{v_1,w_1\}\{v_1,v_2\}\{v_1,v_3\}\{v_1,v_4\} \\	P_2=&\{w_2,v_4\}\{w_2,w_1\}\{w_2,v_3\}\{w_2,w_5\}\{w_2,v_2\}\{w_2,w_4\}\{w_2,v_1\}\{w_2,w_3\} \\
	P_3=&\{w_3,v_4\}\{w_3,w_1\}\{w_3,v_3\}\{w_3,w_6\}\{w_3,v_2\}\{w_3,w_5\}\{w_3,v_1\}\{w_3,w_4\} \\
	P_4=&\{w_4,v_4\}\{w_4,w_1\}\{w_4,v_3\}\{w_4,w_7\}\{w_4,v_2\}\{w_4,w_6\}\{w_4,v_1\}\{w_4,w_5\} \\
	P_5=&\{w_5,v_4\}\{w_5,w_1\}\{w_5,v_3\}\{w_5,w_8\}\{w_5,v_2\}\{w_5,w_7\}\{w_5,v_1\}\{w_5,w_6\} \\
	P_6=&\{w_6,v_4\}\{w_6,w_2\}\{w_6,v_3\}\{w_6,w_1\}\{w_6,v_2\}\{w_6,w_8\}\{w_6,v_1\}\{w_6,w_7\} \\
	P_7=&\{w_7,v_4\}\{w_7,w_3\}\{w_7,v_3\}\{w_7,w_2\}\{w_7,v_2\}\{w_7,w_1\}\{w_7,v_1\}\{w_7,w_8\} \\
	P_8=&\{w_8,v_4\}\{w_8,w_4\}\{w_8,v_3\}\{w_8,w_3\}\{w_8,v_2\}\{w_8,w_2\}\{w_8,v_1\}\{w_8,w_1\} 
	\end{aligned}
	\]
	
	So, concatenating these paths we obtain the Hamiltonian cycle $C_2=P_1P_2\ldots P_8$. 
	
	\item[3)] $\mathbf{1<m<2n.}$
	
	In this case we consider the graph $F_{6,4}$, so in $F_2(F_{6,4})$ we have the following paths:
	\[
	\begin{aligned}
	P'_1=&\{v_4,w_1\}\{v_3,w_1\}\{v_3,v_4\}\{v_2,v_4\}\{v_2,v_3\}\{v_2,w_1\}\{v_1,w_1\}\{v_1,v_2\}\{v_1,v_3\}\{v_1,v_4\} \\	
	P'_2=&\{w_2,v_4\}\{w_2,w_1\}\{w_2,v_3\}\{w_2,w_5\}\{w_2,v_2\}\{w_2,w_4\}\{w_2,v_1\}\{w_2,w_3\} \\
	P'_3=&\{w_3,v_4\}\{w_3,w_1\}\{w_3,v_3\}\{w_3,w_6\}\{w_3,v_2\}\{w_3,w_5\}\{w_3,v_1\}\{w_3,w_4\} \\
	P'_4=&\{w_4,v_4\}\{w_4,w_1\}\{w_4,v_3\}\{w_4,v_2\}\{w_4,w_6\}\{w_4,v_1\}\{w_4,w_5\} \\
	P'_5=&\{w_5,v_4\}\{w_5,w_1\}\{w_5,v_3\}\{w_5,v_2\}\{w_5,v_1\}\{w_5,w_6\} \\
	P'_6=&\{w_6,v_4\}\{w_6,w_2\}\{w_6,v_3\}\{w_6,v_2\}\{w_6,v_1\}\{w_6,w_1\} 
	\end{aligned}
	\]
	
	Therefore, we can concatenate these paths as $C'=P'_1P'_2\ldots P'_6$ to obtain a Hamiltonian cycle in $F_2(F_{6,4})$. 
\end{itemize} 

\noindent\textbf{Complete double vertex graph of $\mathbf{F_{m,n}}$}

As before, we follow the proof of Theorem~\ref{thm:complete_double_vertex}, so we first consider 
the case $m=1$, then the case $m=2\,(n-1)$ and finally the case $1<m<2\,(n-1)$. 
\begin{itemize}
	\item[1)] $\mathbf{m=1.}$
	
	Here we consider the fan graph $F_{1,4}$, so the corresponding paths in $M_2(F_{1,4})$ are:
    \[
	\begin{aligned}
	T'_1=&\{v_1,w_1\}\{v_1,v_1\}\{v_1,v_2\}\{v_1,v_3\}\{v_1,v_4\} \\
	T'_2=&\{v_2,w_1\}\{v_2,v_2\}\{v_2,v_3\}\{v_2,v_4\}\\
	T'_3=&\{v_3,w_1\}\{v_3,v_3\}\{v_3,v_4\}\\
	T'_4=&\{v_4,w_1\}\{v_4,v_4\}\\
	\end{aligned}
	\]
	
	Hence, the concatenation 
	\[
	\begin{aligned}
	C=&\underbrace{\{v_1,w_1\}\{v_1,v_1\}\{v_1,v_2\}\{v_1,v_3\}\{v_1,v_4\}}_{T_1}
	\underbrace{\{v_2,v_4\}\{v_2,v_3\}\{v_2,v_2\}\{v_2,w_1\}}_{\overleftarrow{T_2}} \\
	&\underbrace{\{v_3,w_1\}\{v_3,v_3\}\{v_3,v_4\}}_{T_3}
	\underbrace{\{v_4,v_4\}\{v_4,w_1\}}_{\overleftarrow{T_4}}\{w_1,w_1\} 
	\end{aligned}\]
	is our desired Hamiltonian cycle in  $M_2(F_{1,4})$.

	\item[2)] $\mathbf{m=2\,(n-1).}$
	
	For the graph $F_{6,4}$, we have the following paths in $M_2(F_{6,4})$: 
	\[
	\begin{aligned}[]
	P_1=&\{v_4,w_1\}\{w_1,w_1\}\{v_1,w_1\}\{v_1,v_1\}\{v_1,v_2\}\{v_1,v_3\}\{v_1,v_4\}
	\{v_2,v_4\}\{v_2,v_3\}\{v_2,v_2\}\{v_2,w_1\}\\
	&\{v_3,w_1\}\{v_3,v_3\}\{v_3,v_4\}\{v_4,v_4\} \\
    P_2=&\{w_2,v_4\}\{w_2,w_2\}\{w_2,v_3\}\{w_2,w_1\}\{w_2,v_2\}\{w_2,w_4\}\{w_2,v_1\}\{w_2,w_3\} \\
	P_3=&\{w_3,v_4\}\{w_3,w_3\}\{w_3,v_3\}\{w_3,w_1\}\{w_3,v_2\}\{w_3,w_5\}\{w_3,v_1\}\{w_3,w_4\} \\
	P_4=&\{w_4,v_4\}\{w_4,w_4\}\{w_4,v_3\}\{w_4,w_1\}\{w_4,v_2\}\{w_4,w_6\}\{w_4,v_1\}\{w_4,w_5\} \\
	P_5=&\{w_5,v_4\}\{w_5,w_5\}\{w_5,v_3\}\{w_5,w_2\}\{w_5,v_2\}\{w_5,w_1\}\{w_5,v_1\}\{w_5,w_6\} \\
	P_6=&\{w_6,v_4\}\{w_6,w_6\}\{w_6,v_3\}\{w_6,w_3\}\{w_6,v_2\}\{w_6,w_2\}\{w_6,v_1\}\{w_6,w_1\} 
	\end{aligned}
	\]
	Concatenating these paths we obtain the Hamiltonian cycle $P=P_1P_2\ldots P_6$ in $M_2(F_{6,4})$.

	\item[3)] $\mathbf{1<m<2\,(n-1).}$
	
	Here we consider the graph $F_{5,4}$, and so we have the following paths in $M_2(F_{5,4})$: 
	\[
	\begin{aligned}[]
	P'_1=&\{v_4,w_1\}\{w_1,w_1\}\{v_1,w_1\}\{v_1,v_1\}\{v_1,v_2\}\{v_1,v_3\}\{v_1,v_4\}
	\{v_2,v_4\}\{v_2,v_3\}\{v_2,v_2\}\{v_2,w_1\}\\
	&\{v_3,w_1\}\{v_3,v_3\}\{v_3,v_4\}\{v_4,v_4\} \\
	P'_2=&\{w_2,v_4\}\{w_2,w_2\}\{w_2,v_3\}\{w_2,w_1\}\{w_2,v_2\}\{w_2,w_4\}\{w_2,v_1\}\{w_2,w_3\} \\
	P'_3=&\{w_3,v_4\}\{w_3,w_3\}\{w_3,v_3\}\{w_3,w_1\}\{w_3,v_2\}\{w_3,w_5\}\{w_3,v_1\}\{w_3,w_4\} \\
	P'_4=&\{w_4,v_4\}\{w_4,w_4\}\{w_4,v_3\}\{w_4,w_1\}\{w_4,v_2\}\{w_4,v_1\}\{w_4,w_5\} \\
	P'_5=&\{w_5,v_4\}\{w_5,w_5\}\{w_5,v_3\}\{w_5,w_2\}\{w_5,v_2\}\{w_5,v_1\}\{w_5,w_1\} 
	\end{aligned}
	\]
	Then, the concatenation  $P'=P'_1P'_2\ldots P'_5$ is our desired Hamiltonian cycle in $M_2(F_{5,4})$. 
\end{itemize}

\section{Open problems}
\label{sec:open}

In this paper we have discussed the Hamiltonicity of the double vertex graph and 
the complete double vertex graph of the join graph $G = G_1 + G_2$, where 
$G_1$ and $G_2$ are of order $m\geq 1$ and $n\geq 2$, respectively, and $G_2$ has a Hamiltonian path. 
So, a natural problem is to try to extend these results for $F_k(G)$ and $M_k(G)$. 

\begin{prob}
	Let $G_1$ and $G_2$ be two graphs of order $m\geq 1$ and $n\geq 2$, respectively, and let $G= G_1+G_2$. To study the Hamiltonicity of $F_k(G)$ and $M_k(G)$ for $2<k< n-2$. 
\end{prob}

Also, it can be considered other operations of graphs, such as graph union or graph intersection, 
and some product of graphs.  

\begin{prob}
	Let $G_1$ and $G_2$ be two connected graphs and let $2\leq k\leq n-2$. To study the Hamiltonicity of 
	the $k$-token graph and the $k$-multiset graph of the Cartesian product $G_1\square G_2$, 
	the direct product $G_1\times G_2$ and the strong product $G_1\boxtimes G_2$.
\end{prob}



\end{document}